\documentclass[12pt]{amsart}
\usepackage{amsmath,amssymb, amsthm, graphicx}

\newtheorem{thm}{Theorem}[section]

\newtheorem{prop}[thm]{Proposition}
\theoremstyle{definition}

\theoremstyle{remark}
\newtheorem{rem}[thm]{Remark}

\numberwithin{equation}{section}

\def\beq{\begin{equation} }
\def\eeq{\end{equation}}

\def\cA{{\mathcal A}}

\def\cB{{\mathcal B}}
\def\cD{{\mathcal D}}

\def\cG{\mathcal G}

\def\cF{\mathcal{F}}

\def\ff{\varphi}

\def\ee{\varepsilon}

\newcommand{\NN}{ {\mathbb N} }
\newcommand{\RR}{ {\mathbb R} }
\newcommand{\CC}{{\mathbb C}}

\def\tr{\text{tr}}

\renewcommand{\Im}{\mathrm{Im}\,}
\renewcommand{\Re}{\mathrm{Re}\,}

\begin{document}

\title[Operator-valued semicircular elements]{Operator-valued semicircular elements:
Solving a quadratic matrix
equation with positivity constraints}

\author[W. Helton]{J. William Helton $^{(*)}$}
\address{Math Dept, Univ. of California San Diego, La Jolla Cal, 92093, USA}
\email{helton@ucsd.edu}
\thanks{$^*$ Research supported by the NSF and the Ford Motor Co.}

\author[R. Rashidi Far]{Reza Rashidi Far}
\address{Queen's University, Department of Mathematics and Statistics,
Jeffery Hall, Kingston, ON, K7L 3N6, Canada} \email{reza@mast.queensu.ca}

\author[R. Speicher]{\hbox{Roland Speicher
$^{(\dagger)}$}}
\thanks{$^\dagger\,$Research supported by Discovery and LSI grants from NSERC (Canada) and by
a Killam Fellowship from the Canada Council for the Arts}
\address{Queen's University, Department of Mathematics and Statistics,
Jeffery Hall, Kingston, ON, K7L 3N6, Canada}
\email{speicher@mast.queensu.ca}

\begin{abstract}
We show that the quadratic matrix equation \linebreak $VW + \eta (W)W = I$, for given $V$
with positive real part and given analytic mapping $\eta$ with some positivity preserving
properties, has exactly one solution $W$ with positive real part. Also we provide and
compare  numerical algorithms based on the iteration underlying our proofs.

This work bears on
operator-valued free probability theory,
in particular on the determination of the
asymptotic eigenvalue distribution of band or block random matrices.

\end{abstract}
\thanks{This project got started at the ARCC workshop "Free Analysis" in June 2006. W. Helton and
R. Speicher thank the American Institute of Mathematics for providing a very inspiring
atmosphere. R. Rashidi Far and R. Speicher thank Joachim Cuntz, Siegfried Echterhoff and
the SFB 478 of the Mathematics Institute of the Westf\"alische Wilhelms-Universit\"at
M\"unster for their generous hospitality during their visit in 2006/2007, when main parts
of this research were conducted.}

\maketitle

\section{Introduction}

This paper does two things. One concerns free probability and random matrix theory. The
other concerns iterative methods for solving matrix equations of the form
 \beq \label{eq:VW} VW + \eta
(W)W = I, \eeq where $V$ is a given matrix with positive real part, $\eta$ is a
'positivity preserving' linear or analytic mapping, and we are looking for a solution $W$
with positive real part.

The second topic bears directly on the first. Our main motivation for considering this
kind of equation comes from free probability theory. We will now briefly review the
relevance of this equation in the free probabilistic context of operator-valued
semicircular elements. This will only serve as a motivation and is not essential for our
main statements about solving Eq. \eqref{eq:VW} in section 2. For a general introduction
to free probability theory, see, e.g., \cite{NS}.

Operator-valued semicircular elements \cite{Voi,Spe} play an important role in free
probability
and random matrix theory. In particular,
a big class of random matrices
(band matrices, block matrices) can asymptotically be described by such elements.
This fundamental observation was made by Shlyakhtenko
in \cite{Shl}; extensions and the relevance of this from the point of view
of electrical engineering problems were treated in \cite{ROBS}.

The main problem in a random matrix context is the determination of
the asypmptotic eigenvalue distribution. Let us denote by $H(z)$ the
Cauchy transform of the asymptotic eigenvalue distribution $\mu$, given by
$$H(z):=\int_\RR \frac 1{z-t}d\mu(t).$$
This is an analytic function in the upper complex half plane
$$\CC^+:=\{z\in\CC\mid \Im z> 0\}$$
and allows to recover $\mu$ via the Stieltjes inversion formula
$$d\mu(t)=-\frac 1\pi \lim_{\ee\to 0}\Im H(t+i\ee).$$
The theory of operator-valued
semicircular elements tells us that one can get this Cauchy transform as
$$H(z)=\ff(G(z)),$$
where $G(z)$ is an operator-valued function, i.e. $G(z)\in\cA$ for some operator algebra
$\cA$ (with involution $*$, containing the unit $I$), and where $\ff:\cA\to\CC$ is a
given state on this algebra. (The data $\cA$ and $\ff$ are determined by the form of the
considered matrices; e.g., for block matrices as in \cite{ROBS}, $\cA=M_d(\CC)$ are just
the $d\times d$-matrices, for some fixed $d$, and $\ff=\tr_d$ is the normalized trace on
those matrices.) The operator-valued Cauchy transform $G(z)$ is determined by the
equation \beq \label{eq:G} z G(z)=I+\eta(G(z))G(z)\qquad (z\in \CC^+), \eeq together with
its asymptotic behaviour $G(z)\sim \frac 1z I$ for $z\to\infty$.
%
%$I\in\cA$ is here the unit in
Here we are given the completely positive linear map $\eta:\cA\to\cA$; it contains the
information about the covariance of the considered random matrix or operator-valued
semicircular element.

One
of the main problems, both from a conceptual
and a numerical point of view,
is that, for a fixed $z\in\CC^+$,
the equation \eqref{eq:G} has not just one, but
many solutions. To isolate the correct root is not obvious at all.

We shall use the  common notation:
the real and imaginary parts of an operator $W$ are defined as
$$\Re W:=\frac 12 (W+W^*)  \qquad \Im W:=\frac{ 1}{2i} (W-W^*).$$
Since
$G$ is also the operator-valued Cauchy transform
of a semicircular element $s$,
$$G(z)=E[\frac 1{z-s}]$$
(for some $C^*$-algebra $\cB$ with $s\in\cB$, $s=s^*$,
and a conditional expectation $E:\cB\to\cA$.),
 we have for $\Im z\geq 0$ that
\begin{align*}
\Im G(z)&=\frac 1{2i}(G(z)-G(z)^*)\\
&=\frac 1{2i} E[\frac 1{z-s} -\frac 1{\bar z-s}]\\
&=-\Im z \cdot E[\frac{1}{z-s}\frac{1}{\bar z-s}]\\
&\leq - \frac{\Im z}{\bigl(\vert z\vert+\Vert s\Vert\bigr)^2}\cdot I;
\end{align*}
since
$$\frac{1}{z-s}\frac{1}{\bar z-s}
\geq \frac{1}{\bigl(\vert z\vert+\Vert s\Vert\bigr)^2}
%\cdot I
$$
and a conditional expectation is positive.

This means that we are looking, for any $z\in\CC^+$,
for a solution $G(z)$ of \eqref{eq:G}
which has the property that its imaginary part is negative
and stays away from zero (which is the same as saying that
its imaginary part is negative and invertible).
Here we will show that actually
there exists exactly one solution of \eqref{eq:G} with this property,
and that one can get
it by an appropriately chosen iteration procedure.

If we write $G(z)=-i W(z)$, then $\Im G(z)\leq -\ee I$ gets replaced by the nicer
condition $\Re W(z)\geq \ee I$. Thus we will work with $W(z)$ instead of $G(z)$ in the
following. The equation \eqref{eq:G} reads in terms of $W$ as
\beq -izW(z) + \eta
(W(z))W(z) = I \eeq

As it turns out the linearity of $\eta$ is not crucial (as long as it is analytic, and
has certain positivity and boundedness properties). Thus we will treat the problem
directly in this more general setting. Furthermore, one can also replace the complex
number $z\in\CC^+$ by an arbitrary element $Z\in\cA$ with positive imaginary part. Again
we replace this by the condition of positive real part by going over to $V=-iZ$, ending
up with the Eq. \eqref{eq:VW}.

\section{Main Theorem}

\subsection{Setting of our problem}
We work in a $C^*$-algebra, denoted $\cA$.
Interesting examples arise already in the
finite-dimensional case,
thus the reader might just take $\cA$ as a matrix algebra of
complex $d\times d$-matrices, $\cA=M_{d}(\CC)$ for some $d\in\NN$.

For a selfadjoint operator $A$ we mean with $A \geq 0$ that its spectrum is contained in
$[0,\infty)$. Let $\cA_+$ denote the strict right half plane of $\cA$,
$$\cA_+:=\{W\in\cA\mid \text{$\Re W\geq \ee I$ for some $\ee>0$ }\}.$$
Note that $W\in\cA_+$ if and only if $\Re W\geq 0$ and $\Re W$ is invertible.

%MOVED DEF OF Re FORWARD

Furthermore, we are given an analytic mapping
$$\eta:\cA_+\to\cA_+,$$
which is bounded on bounded domains in $\cA_+$, i.e.
$$\sup\{\Vert \eta(W)\Vert\mid W\in\cA_+,\,\Vert W\Vert\leq r\}<\infty$$
for any $r>0$.

For a given $V\in\cA_+$ we consider our key equation: \beq \label{eq:Weq} VW + \eta (W)W
= I. \eeq We are looking for a solution $W\in\cA_+$.

Note that the solution of this equation is the same as a fixed point of the map $W\mapsto
\cF_V(W)$, where
\begin{equation}\label{eq:cF}
\mathcal{F}_V (W) := [ V + \eta (W)]^{-1}.
\end{equation}

Here is our main result.

\begin{thm}
\label{thm:converge} For fixed $V\in\cA_+$, there exists exactly one solution $W\in\cA_+$
to \eqref{eq:Weq}; this $W$ is the limit of iterates
 $$W_n = \mathcal{F}_V^n (W_0)$$
 for any $W_0\in\cA_+$.
Furthermore, we have that
$$\Vert W\Vert\leq \Vert(\Re V)^{-1}\Vert$$
and
$$\Re W\geq \frac{1}{m^2\cdot\Vert (\Re V)^{-1}\Vert}I,$$
where
$$m:=\Vert V\Vert +
\sup\{\Vert \eta(W)\Vert\mid W\in\cA_+,\,\Vert W\Vert\leq \Vert(\Re V)^{-1}\Vert\}.
$$
\end{thm}

\begin{rem}
The uniqueness part of our theorem can be used to give an easy direct proof of Prop. 5.6 in
\cite{HT}. There two operator-valued Cauchy transforms $G(\lambda)$ and $G_n(\lambda)$ are
considered (where $\lambda$ might also be a matrix, with positive imaginary part), and it is
shown that both $G_n(\lambda)$ and $G((\Lambda_n(\lambda))$ fulfill the same equation, which
is, after our rescaling, of the form \eqref{eq:Weq}. Since both $G(\Lambda_n(\lambda))$ and
$G_n(\lambda)$ satisfy (as Cauchy transforms at some value of the argument) the right
positivity condition, they both must agree with the unique solution, given by our theorem, and
thus $G(\Lambda_n(\lambda))=G_n(\lambda)$.
\end{rem}

\subsection{Related Topics}
\def\cS{\mathcal S}

Suppose $\eta$ maps selfadjoint elements to selfadjoint elements;
then it maps positive elements $\cS\cA_+$  of $\cA$ to
 $\cS\cA_+$.
 Suppose further that $V >0$.
Then, if $\eta$ is linear,
the map $\cF_V$ restricted to $\cS\cA_+$ is a
monotone decreasing map, that is,  if $W_1 \geq W_2 \geq 0$,
then $\cF_{V}(W_1) \leq    \cF_{V}(W_2)$.
There are results in \cite{RR02} and \cite{HS05}
yielding fixed points in this  case
 with proofs quite different than here.
A list of applications and existing results on special cases
is in these papers.

\section{Contraction Maps and Proofs}
We will prove our theorem by applying Banach's fixed point theorem to the map $\cF_V$.
One should note that $\cF_V$ is usually not a contraction in the given operator norm on
$\cA_+$. Analyticity of our mapping, however, provides us with another metric in which we
have the contraction property.

On $\cA_+$, since it is a domain in a Banach space, is the well known Carath\'eodory
metric, one of the biholomorphically invariant metrics on $\cA_+$. See for \cite{Ha03}
for an excellent exposition of this and material germain to our treatment here. The
crucial point is that strict holomorphic mappings on such domains are automatically
strict contractions in this metric, and thus Banach's fixed point theorem guarantees a unique
fixed point of such mappings. For the reader's convenience we recall here the relevant
theorem, due to Earle and Hamilton \cite{EH}.

\begin{thm}{\bf (Earle-Hamilton)}
Let $\cD$ be a nonempty domain in a complex Banach space $X$ and let $h:\cD\to\cD$ be a
bounded holomorphic function. If $h(\cD)$ lies strictly inside $\cD$ (i.e., there is some
$\epsilon>0$ such that $B_\epsilon(h(x))\subset\cD$, whenever $x\in\cD$, where
$B_\epsilon(y)$ is the ball of radius $\epsilon$ about $y$), then $h$ is a strict
contraction in the Carath\'eodory-Riffen-Finsler metric $\rho$, and thus has a unique
fixed point in $\cD$. Furthermore, one has for all $x,y\in\cD$ that $\rho(x,y)\geq m
\Vert x-y\Vert$ for some constant $m>0$, and thus $(h^n(x_0))_{n\in\NN}$ converges in
norm, for any $x_0\in\cD$, to this fixed point.
\end{thm}

We will now apply this to our situation. The main point will be to check that $\cF_V$ is
well-defined and maps suitably chosen subsets $R_b\subset\cA_+$ strictly into itself.
(Note that we do not claim that $\cF_V$ is a strict contraction on $\cA_+$ itself.)

For $b>0$, let us define
$$R_b:=\{W\in \cA_+\mid \Vert W\Vert<b \}\subset \cA_+$$

\begin{prop}
\label{prop:main}
\begin{enumerate}
\item %1
Any fixed point $W$ of $\cF_V$ satisfies Eq. \eqref{eq:Weq}.

\item %2
If $\,\Re \eta (W)\geq 0$ and  $V\in\cA_+$, then $V + \eta (W)$ is invertible in $\cA$,
thus $\cF_V(W)\in\cA$ is well-defined, and $\Re \mathcal{F}_V (W) \geq 0$. Furthermore,
we have
$$\Vert \cF_V(W)\Vert\leq \Vert{(\Re V})^{-1}\Vert.$$

\item %5
\label{it:strict} For $V\in\cA_+$ and $b > \Vert(\Re V)^{-1}\Vert > 0$, the map
$\mathcal{F}_V$ is a bounded holomorphic map on $R_b$ and takes $R_b$ strictly into its
interior. We have for $W\in R_b$ that
$$\Re \cF_V(W)\geq \frac{1}{{m_b^2}\cdot\Vert(\Re V)^{-1}\Vert}I$$
with
$$m_b:=\Vert V\Vert \; + \
\sup\{\Vert \eta(W)\Vert\mid W\in\cA_+,\,\Vert W\Vert\leq b\}.$$

\item

$\cF_V$ maps, for any $V\in\cA_+$, the strict upper half plane of $\cA$ into itself,
$$\cF_V:\cA_+\to\cA_+.$$
It is a bounded holomorphic map there.

\end{enumerate}
\end{prop}

\begin{proof}
(1) This is obvious.

(2) Let $A$ and $B$ be the real and imaginary part of $\eta(W)$, respectively, i.e.,
$\eta(W)=A + iB$ with $A=A^*$ and $B=B^*$. By our assumption, $A\geq 0$. Put $K:=V + \eta
(W)$. Then we have:
$$\Re K=\Re V +A \geq \Re V.$$
Since $V\in\cA_+$, the real part of $K$ is bounded away from zero by a multiple of
identity, and thus $K$ is invertible as a bounded operator. Furthermore, the norm of the
inverse $K^{-1}=\cF_V(W)$ can be bounded by $\Vert (\Re V)^{-1}\Vert$. (For more details
about these statements, see Lemma 3.1 in \cite{HT}.)

To see the positivity of the real part of the inverse, we calculate:
\begin{align*}
 2\cdot \Re \cF_V (W) &= [V + A + iB]^{-1}
+ [V + A + iB ]^{\ast-1} \\
&= [V + A + iB]^{-1} (2 A + V+V^*)
[V^* + A - iB]^{-1}\\
&= 2 \mathcal{F}_V(W)(A + \Re V) \mathcal{F}_V (W)^*\\
&\geq 2\mathcal{F}_V(W)\cdot \Re V\cdot \mathcal{F}_V (W)^*\\
&\geq 0.
\end{align*}

(3) The map $W\mapsto K(W):=V + \eta (W)$ is analytic by the analyticity of $\eta$, and
$\cF_V(W)=K(W)^{-1}$ exists -- and is thus also analytic -- for $W\in\cA_+$, by (2) (note
that the fact that $\eta$ preserves the positivity of the real part implies that the
assumptions of (2) are satisfied). Furthermore, by the norm estimate from (2), we have
for all $W\in R_b$ that \beq \label{eq:norm-estimate} \Vert \cF_V(W)\Vert\leq \Vert (\Re
V)^{-1}\Vert< b. \eeq Thus $\cF_V$ is a bounded holomorphic map on $R_b$, with image in
$R_b$. To see that the image lies strictly in $R_b$, we have to see that $\cF_V(W)$ stays
away from the boundary of $R_b$ by some $\ee$-amount. By \eqref{eq:norm-estimate} we stay
away from the boundary $\Vert W\Vert=b$ by at least $(b-\Vert(\Re V)^{-1}\Vert)$. It
remains to consider the part of the boundary described by $\Re W=0$. By refining the last
inequality of the calculation from (2) in our present setting we have
\begin{align*}
\Re \mathcal{F}_V(W) & \geq
\mathcal{F}_V(W)\cdot \Re V\cdot \mathcal{F}_V(W)^\ast\\& =
[K(W)^\ast\cdot(\Re V)^{-1}\cdot K(W)]^{-1} \\
& \geq \frac{1}{\lVert K(W)^* \cdot(\Re V)^{-1}\cdot K(W) \rVert} I \\&\geq
\frac{1}{{m_b^2}\cdot\Vert(\Re V)^{-1}\Vert}I
\end{align*}
since $\lVert K(W) \rVert \leq m_b$ for all $W \in R_b$. Thus we stay away from the
boundary $\Re W=0$ by at least $1/({m_b^2}\cdot\Vert(\Re V)^{-1}\Vert)$.

(4) This follows from the fact that
$$\cA_+=\bigcup_{b>0}R_b;$$
note that the estimate $\Vert \cF_V(W)\Vert\leq \Vert(\Re V)^{-1}\Vert$ does not depend
on $b$.
\end{proof}

\begin{proof}[Proof of Theorem \ref{thm:converge}]
By the Earle-Hamilton Theorem, each $R_b$ with $b>\Vert(\Re V)^{-1}\Vert$ contains
exactly one fixed point of \eqref{eq:Weq}. The estimates for the norm and the real part
of $W$ follow from the corresponding estimates in parts (2) and (3) of Prop.
\ref{prop:main}.
\end{proof}

\begin{rem}
1) Since the application of Stieltjes inversion formula asks for $z$ very close to the
real axis, one might be tempted to try to solve \eqref{eq:G} directly for real $z$. Of
course, most of the above statements then break  down.
In particular, one should consider
the map $\cF_V$ for $V$ with $\Re V\geq 0$ on the domain $\Re W>0$ instead of $\cA_+$. In
the infinite dimensional case, those two notions are not the same, and using $\Re W>0$
presents problems. ($\cF_V$ is not even well-defined there in general, for $\Re V\geq
0$.) In the finite dimensional case, however, $\Re W>0$ just says that all eigenvalues of
$\Re W$ are positive and different from zero, thus $\Re W$ is positive and invertible,
thus $W\in\cA_+$ and one can extend some of the above reasoning to the set
$\{W\in\cA\mid \Re W>0\}$. This domain is mapped, under $\cF_V$ into itself, but now of
course not necessarily strictly. This implies (see Proposition 6 in \cite{Kran} ) that,
we still get a contraction in the Caratheodory metric, but the contraction constant is
not necessarily less than 1. We now state this with formulas.
 Let
$\cA$ be finite dimensional (i.e., some matrix algebra).
Let $\rho$ denote the
Caratheodory metric on the set $\{W\in\cA\mid \Re W>0\}$.
 Then for $\Re V \geq 0$, there
is $c_V \leq 1$ such that
$$\rho (\mathcal{F}_V (W), \mathcal{F}_V (\widetilde{W}))
\leq c_V \; \rho( W, \widetilde{W})$$ if $\Re W > 0 \text{ and } \Re \widetilde{W} > 0$.

2) One can generalize our considerations from \eqref{eq:G} to the equation
\begin{equation}\label{eq:general}
G(z)=G_1[z-\eta(G(z))].
\end{equation}
The latter describes \cite{Spe,Voi} the operator-valued Cauchy transform of the sum $x+s$
where $s\in\cB$ is an operator-valued semicircular element as before (with covariance
function $\eta$) and $x=x^*\in\cB$ is an element which is free from $s$ with respect to
the conditional expectation $E:\cB\to\cA$;
$$G_1(z):=E[\frac 1{z-x}]$$
is the given operator-valued Cauchy transform of $x$. Note that for $x=0$ we have
$G_1(z)=1/z$ and \eqref{eq:general} reduces to the fixed point version of \eqref{eq:G}.
In the scalar-valued case $\cA=\CC$, equation \eqref{eq:general} was derived by Pastur
\cite{Pas}, describing a 'deformed semicircle'.

The same arguments as before show that for any fixed $z$ with positive imaginary part
there exists exactly one solution of \eqref{eq:general} whose imaginary part is strictly
negative. This unique solution can be obtained as the limit of iterates of the mapping
$G\mapsto G_1[z-\eta(G)]$, for any initial $G$ with strictly negative imaginary part.

\end{rem}

\section{Numerical considerations}

\subsection{Iteration and averaging}
Let us come back to our motivating equation \eqref{eq:G} for the operator-valued Cauchy
transform $G(z)$ of an operator-valued semicircular element. According to Theorem
\ref{thm:converge}, we can get the wanted solution $G=-iW$ of Eq. \eqref{eq:G} by
iterating the mapping $W\mapsto\cF_z(W):=[-izI+\eta(W)]^{-1}$, starting with any
$W_0\in\cA_+$. Even though this seems to completely solve our problem, it turns out that
numerically this might not work too well. Namely, since we want to invoke the Stieltjes
inversion formula to recover the wanted eigenvalue distribution from $G(z)$ we need $z$
very close to the real axis. Then $\cF_z$ is still a contraction
(in the Caratheodory
metric), the contraction constant, however, might be very close to 1, and the convergence
could be extremely slow. In typical numerical examples, there was, for $z$ very close to
the real axis, no numerically observable convergence to a fixed point at all. (In some
cases it looked as if $\cF_z$ would have a limiting 2-cycle.) To improve on this, we
invoked also an averaging procedure along with our iteration. In the simplest case, we
replaced the iteration algorithm
$$W\mapsto \cF_z(W)$$
by
$$W\mapsto \cG_z(W):=\frac 12 W+\frac 12\cF_z(W).$$

Note that $\cG_z$ has essentially the same properties  as $\cF_z$; in particular, the fixed
points are the same, and $\cG_z$ is a bounded holomorphic map on $R_b$ and takes $R_b$
strictly into its interior, as follows from the following lemma.

\begin{prop}
Let $\cD$ be a convex and bounded domain in a
Banach space $X$ and $\cF:\cD\to\cD$ a bounded holomorphic
function. Assume that $\cF(\cD)$ lies strictly inside $\cD$.
Put
$$\cG(W):=\frac 12 W+\frac{1}{2}\cF(W).$$
Then $\cG:\cD\to\cD$ is bounded and holomorphic and
maps $D$ strictly into itself.
Thus $\cG$ iterations applied to any $W^0$ inside $\cD$
  converge to a fixed point for both $\cG$ and $\cF$.
\end{prop}

\begin{proof}
Only the fact that $\cG(\cD)$ lies strictly inside $\cD$ is not directly clear.
Let $\ee>0$ be such that $B_\ee(\cF(W))\subset\cD$ for all $W\in\cD$.
We claim that
$B_{\ee/2}(\cG(W))\subset \cD$ for all $W\in \cD$.
To see this assume there is a $W\in\cD$ and
an $r\in X$
with $\Vert r\Vert\leq \ee/2$ such that $\cG(W)+r\not\in\cD$.
But we know that $\cF(W)+2r\in\cD$ and thus, since $\cD$ is convex, also
$$\cG(W)+r=\frac12 (\cF(W)+2r)+\frac 12 W\in\cD.$$

The convergence of $\cG$ iterates,
follows immediately from the Earle-Hamilton Theorem.
\end{proof}

%So we can apply the Earle-Hamilton Theorem to $\cG_z$
%and see that
As stated iterations of
$\cG_z$ must converge to
the same fixed point as iterations of $\cF_z$.
Numerically, in all our considered reasonable examples,
the speed of convergence for the averaged iteration was
 substantially faster than for
the plain iterations and allowed a fast
 numerical determination of the
desired fixed point.

Let us, however, point out that
there is no theoretical reason that the averaged iteration has a better
contraction constant than the plain iteration (actually it might even be worse), and that
one can produce artificial examples where the averaging does not improve the plain
iteration. For example, in the case of $2\times 2$-matrices where
$$\eta(W)= AWA,\qquad\text{with} \qquad A=\begin{pmatrix} 0&1\\1&0\end{pmatrix},$$
and for initial condition
$$W_0=\begin{pmatrix} 2&0\\0&1/2
      \end{pmatrix},$$
there is no clear improvement of the averaged iteration over the plain iteration
for $z$ close to 0.

Of course, the usual algorithm of choice for solving equations like \eqref{eq:G} or
\eqref{eq:Weq} would be a kind of Newton's method
applied at later iterations.
This has much faster convergence
properties (second order local convergence compared to first order convergence rate of
our iteration). However, Newton's method does in general not respect the positivity
requirement for our solution. So even if we start in $\cA_+$, it might happen that the
Newton's method leads out of $\cA_+$ and converges to another root which is not in
$\cA_+$. One might expect that iterating (or iterating and averaging) for a while will
get in the basin of attraction of the correct point and then switching to the Newton's
method should be fairly safe and fast.

Whereas we do not have any theoretical results to justify averaging or a modification of
the Newton's method in all cases, we want to point out that the uniqueness part in our Theorem allows
us to check in a concrete example whether our algorithm of choice works or not. We only
have to check whether the solution produced by our algorithm is in $\cA_+$.

\subsection{Numerical Examples}
We consider the block Toeplitz random matrix
\begin{eqnarray}\label{eq:3x3mtx}
X=\frac{1}{\sqrt{3N}}
 \begin{bmatrix} A &B &C \\
B&A&B\\ C&B&A
\end{bmatrix},
\end{eqnarray}
where $A,B,C$ are independent selfadjoint $N\times N$ matrices with i.i.d. entries of unit
variance. In \cite{ROBS} it was shown that, in the limit $N\to\infty$, the Cauchy-transform $H$
of the eigenvalue distribution of $X$ is given by $H(z)=\tr_3(G(z))$ where $\tr_3$ denotes the
normalized trace on $3\times 3$-matrices and where the operator-valued Cauchy transform $G(z)$
is a $3\times 3$-matrix of the form
\begin{eqnarray}\label{eq:G3x3} G=\begin{bmatrix}
   f & 0 & h  \cr 0 & g & 0  \cr
     h & 0 & f \cr  \end{bmatrix}.
\end{eqnarray}
$G$ satisfies Eq. \eqref{eq:G}, where $\eta$ acts as follows:
\begin{eqnarray}\label{eq:eta3x3}
 \;\eta(G)=\frac13\begin{bmatrix}  2\,f + g & 0 & g + 2\,h \cr 0 & 2\,f +
   g + 2\,h & 0 \cr g + 2\,h & 0 & 2\,f + g   \end{bmatrix}.
\end{eqnarray}

In order to get the eigenvalue distribution of $X$ we have to solve Eq. \eqref{eq:G} for
$z$ running along the real axis. We solve \eqref{eq:G} either by our modified iteration
procedure or by Newton's algorithm. In both cases we start from the center where $z_0$ is
an imaginary number close to the origin ($z_0=10^{-9}\cdot i$ for the following results)
and we use the solution at each $z_n=z_0+d\cdot \Delta\cdot n$ as the initial point for
the next $z_{n+1}=z_0+d\cdot \Delta\cdot \left(n+1\right)$ where $\Delta \in \RR$ is the
resolution step size in calculating the spectrum and $d=\mp 1$ depending on negative or
positive side of the spectrum. In other words, we calculate the spectrum in two phases
both starting from the center, one for the positive side of the spectrum ($d=1$), then
for the negative side ($d=-1$).

\subsubsection{Solution by iteration}

\begin{figure}\centering{
\includegraphics[width=12cm]{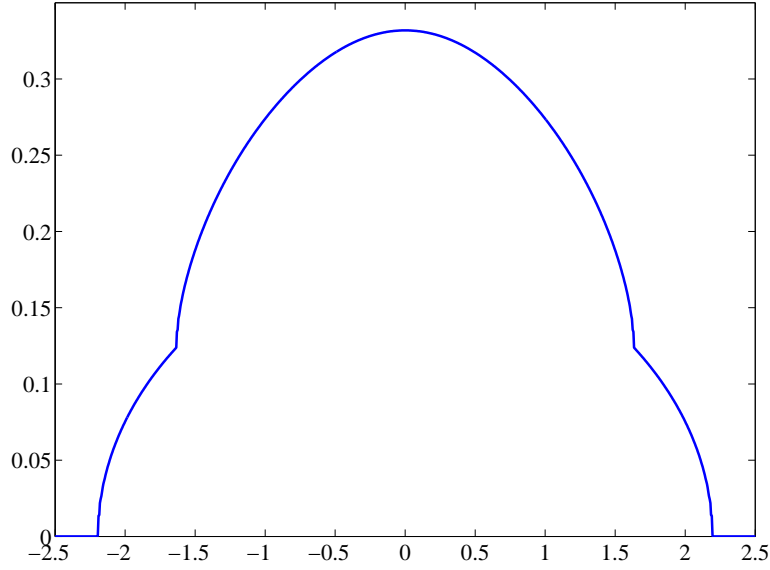}}
\caption{\label{fig:Iteration}Iteration method calculation result for the asymptotic
eigenvalue distribution of the $3\mathbf{x}3$ block Toeplitz matrix as in (\ref{eq:3x3mtx}).}
\end{figure}

With the convention $W(z)=iG(z)$, we use
\begin{eqnarray}\nonumber
W_0(z)=\left[\begin{array}{ccc}
1-0.1\cdot i& 0& 1\\0& 1-0.1\cdot i& 0\\1& 0& 1-0.1\cdot i \end{array}
\right]\cdot i
\end{eqnarray}
as the initializing matrix at center ($z=10^{-9}\cdot i$) for the iteration method.
Figure~\ref{fig:Iteration} depicts the calculated spectrum which matches completely with
the true spectrum \cite{ROBS}.

\subsubsection{Failure of Newton's method}
\begin{figure}\centering{
\includegraphics[width=12cm]{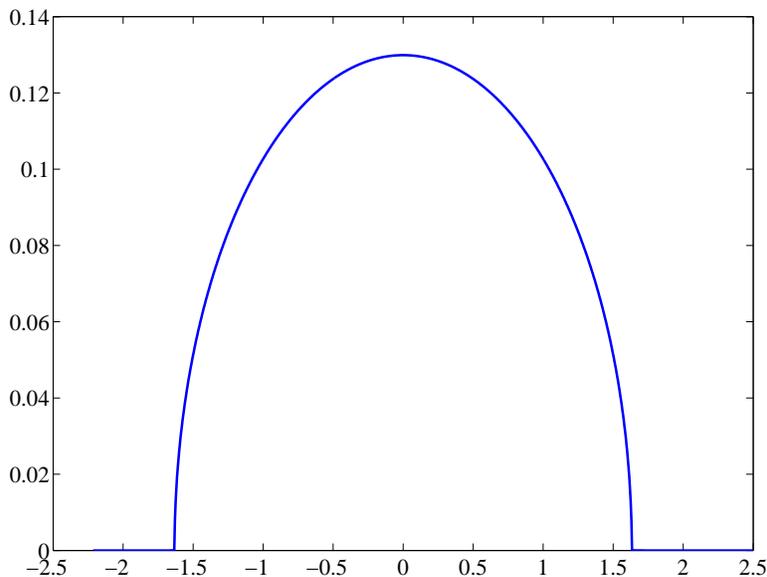}}
\caption{\label{fig:Newton}Newton's method calculation result for the eigenvalue distribution of the $3\mathbf{x}3$ block Toeplitz matrix as in (\ref{eq:3x3mtx}).}
\end{figure}

Replacing (\ref{eq:G3x3}) and (\ref{eq:eta3x3}) in (\ref{eq:G}), one reaches to the following system of equations:
\begin{eqnarray}\label{eq:3x3sys}
    z f &=& 1 + \frac{g\,\left( f + h \right)  +
    2\,\left( f^2 + h^2 \right) }{3},\\
  z g &=& 1 + \frac{g\,
       \left( g + 2\,\left( f + h \right)  \right)
         }{3}, \\ zh  &=&
   \frac{4\,f\,h + g\,\left( f + h \right) }{3}.
\end{eqnarray}
We used Newton's method to solve this system numerically, starting with the initial
values
\begin{eqnarray}\nonumber
\left[
\begin{array}{c}
f_0\\
g_0\\
h_0
\end{array}
\right]=
\left[
\begin{array}{c}
1-0.1\cdot i\\
1-0.1\cdot i\\
1
\end{array}
\right],
\end{eqnarray}
for $z=10^{-9}i$ on the imaginary axis and then using the fixed point for one $z_n$ as
the initial point for the next $z_{n+1}$. The result is shown in Figure~\ref{fig:Newton}.
Clearly this is not the desired result and the algorithm has converged to some different
roots, which in this case do not yield a distribution function (though it is positive
everywhere). In other words, using Newton's method, one has to check the positivity of
the final result to make sure that it has converged to the desired solution. Of course,
in our example the $3\times 3$ matrix $G$ produced by Newton's algorithm is not positive.

\end{document}